\documentclass[11pt,leqno]{article}
\usepackage{amsthm,amsbsy,amsfonts,amssymb,amsmath,amscd}

\newcommand{\Q}{\mathbb Q}
\newcommand{\Z}{\mathbb Z}

\newcommand{\C}{\mathbb C}

\newcommand{\F}{\mathbb F}

\newcommand{\Hom}{{\rm Hom}}

\newcommand{\res}{\rm res}

\newcommand{\Sh}{\rm Sha}
\newcommand{\ScrG}{\mathcal G}
\usepackage{latexsym}

\makeatletter
\def\l@section{\@tocline{1}{4pt}{1pc}{}{}}
\def\l@subsection{\@tocline{2}{0pt}{2pc}{5pc}{}}
\makeatother

\title{Global Galois Symbols on $E \times E$}
\author{Dinakar
Ramakrishnan\footnote{Partially supported by the Simons Foundation (grant \# 523557)}}
\date{}
\begin{document}
\maketitle \overfullrule=5pt
\newcounter{mysubsection}[section]
\setcounter{mysubsection}{0}
\newcommand{\mysubsection}[1]{\refstepcounter{mysubsection} %
\vspace{0.5em} \noindent {\large\bf
\arabic{section}.\arabic{mysubsection} \hspace{0.4em} #1}
\par\vspace*{.2truecm}
}

\begin{flushright}
{\it In memory of Jacob Murre: friend, teacher and coauthor}
\end{flushright}

\addcontentsline{toc}{section}{Introduction}
\section*{Introduction}

\bigskip

Let $E$ be an elliptic curve over a number field $F$, $\overline F$ an
algebraic closure of $F$, and $A$ the abelian surface $E \times E$. Denote by $CH_0(A)$ the Chow group of zero cycles on $A$ and $CH_0(A)^0$ its degree zero part. Denote by $T_F(A)$
the {\it albanese kernel}, i.e., the kernel of the albanese map $CH_0(A)^0 \rightarrow {\rm Alb}(A)$ over $F$, where Alb$(A)$ denotes the albanese variety of $A$. Since $A$ has positive geometric genus, one knows by Mumford and Bloch that over $\C$, $T(A)$ is very large, in fact infinite dimensional in a suitable sense. On the other hand, a deep conjecture of Beilinson (\cite{Be}) and Bloch (\cite{Bl}) asserts that over a number field $F$, $T_F(A)$ is a torsion group. Moreover, a conjecture of Bass on $K_0$ of schemes of finite type over $\Z$ predicts, thanks to Riemann-Roch without denominators (\cite{CT-R}), that this group is finitely generated, so $T_F(A)$ should even be finite. (See \cite{LaS} and \cite{LaR} for positive results on the torsion subgroup for a product of certain classes of elliptic curves.)
The expected finiteness of $T_F(X)$ is unknown for any surface $X$ over $F$ with $p_g >0$, and our aim is more modest. It began with Jaap (Murre) desiring a way to test the conjecture in an example, at least an aspect of it. A simple consequence of the expected finiteness is that for all primes $p$ large enough, $T_F(X)/p$ should be zero. We verify, for $X=A$, that the image of $T_F(A)/p$ in Galois cohomology $H^2(F, E[p] \otimes E[p])$ (see below) is zero for a set of $p$ of density $1$ (resp. $1/2$) if $E$ is non-CM (resp. CM). To be precise we show this for all but a finite number of $p$ where $E$ has good {\it ordinary} reduction. The supersingular case is subtle and we hope to discuss it elsewhere.

Most of the work on this was done jointly wih Jaap Murre. Sadly we did not finish writing it up before his passing. {\it Sorry voor de vertraging, Jaap! Vergeef mij alstublieft mijn neiging tot uitstel!}

It should be noted that this article is a sequel to our earlier joint paper \cite{MR}, which dealt with local Galois symbols on $E \times E$. We will freely use the results and definitions from that paper, which contains some of the key preliminaries, and implicitly the references cited in it.

Let us recall briefly how the Galois cohomology enters.
Denote by $E[p]$ the group of
$p$-division points of $E$ in $E(\overline F)$. To any
$F$-rational point $P$ in $E(F)$ Kummer theory yields a
class $[P]_p$ in the Galois cohomology group $H^1(F, E[p])$,
represented by the $1$-cocycle
$$
\beta_p: {\rm Gal}(\overline F/F) \rightarrow E[p], \, \,
\sigma \mapsto \sigma\left(\frac{P}{p}\right) - \frac{P}{p}.
$$
Here $\frac{P}{p}$ denotes any point in $E(\overline F)$ with
$p\left(\frac{P}{p}\right) = P$. Given a pair $(P,Q)$ of
$F$-rational points, one then has the cup product class
$$
[P,Q]_p : = \, [P]_p \cup [Q]_p \, \in \, H^2(F,
E[p]^{\otimes 2}).
$$

Any such pair $(P,Q)$ also defines a $F$-rational algebraic cycle
on the surface $E \times E$ given by
$$
\langle P,Q \rangle : = \, [(P,Q) - (P,0) - (0,Q) + (0,0)],
$$
where $[ \cdots ]$ denotes the class taken modulo {\it rational
equivalence}. It is clear that this {\it zero cycle of degree
zero} defines, by the parallelogram law, the trivial class in the
Albanese variety $Alb(A)$. So $\langle P,Q\rangle$ lies
in the {\it Albanese kernel} $T_F(A)$. It is a known, but
non-obvious, fact that the association $(P,Q) \to [P,Q]_p$
depends only on $\langle P,Q\rangle$, and thus results in the {\it
Galois symbol map} on the subgroup $ST_F(A)$ generated by the symbols $\langle P,Q \rangle$.
One knows (cf. \cite{MR}) that the entire $T_F(A)$ is generated by norms
$N_{K/F}(\langle P', Q'\rangle)$ for finite extensions $K/F$ with $P', Q' \in E(K)$, which is compatible with the norm map (corestriction) on Galois cohomology. One uses it to define the Galois symbol map
$$
s_p: \, T_F(E \times E)/p \, \rightarrow \, H^2(F,
E[p]^{\otimes 2}).
$$
For the precise technical definition of $s_p$ we refer to Definition 2.5.3 of \cite{MR}. In essence, $s_p$ is the restriction to the Albanese kernel $T_F(E\times E)$ of the cycle map over $F$ with target the continuous cohomology (in the sense of Jannsen) with $\Z/p$-coefficients.

It is reasonable to focus on this map and try to show (for large $p$) the triviality of the image of $T_F(A)/p$, since a conjecture of Somekawa and Kato (\cite{So}) asserts that the Galois symbol map is always injective (\cite{So}). This conjecture is still open, though known locally for $p$ good ordinary (\cite{MR}) for $p$ multiplicative (\cite{Y}); we will make use of these local results below.

\medskip

\noindent{\bf Theorem A} \quad \it Let $E$ be an elliptic curve over a number field $F$. Then for all but a finite set of rational primes $p$ where $E$ has ordinary reduction,
$$
s_p(T_F(A)/p) \, = \, 0.
$$
\rm

\medskip

When $E$ is not of CM type, a theorem of Serre \cite{Serre2}, section 8, asserts that the set of ordinary primes for $E$ has density one.    When $E$ has CM by an order in an imaginary quadratic field $K$, it has good ordinary reduction at almost all primes $p$ which split in $K$.

There is a $\Q_p$-version of this Theorem, for $E$ of CM type, in \cite{LaR}, Thm. 0.2.

\medskip

The next result goes in the other direction:

\medskip

\noindent{\bf Theorem B} \quad \it Let $E$ be an elliptic curve over a number field $F$. Let $p$ be any prime where $E$ has good ordinary reduction. Then there exists a finite extension $K$ of $F$ (depending on $p$) such that
$$
T_K(A)/p \, \ne \, 0.
$$
\rm

\medskip

Some of the classes we construct in $T_K(A)$ are torsion, but the others are more mysterious.

\medskip

Note that there is a Galois module decomposition
$$
E[p]^{\otimes 2} \, \simeq \, {\rm sym}^2(E[p]) \oplus \Lambda^2(E[p]),
$$
where the second module on the right is by the Weil pairing the one dimensional representation $\mu_p$
(cyclotomic character mod $p$). It is known (see \cite{MR}, for example) that the composite map $\pi_p\circ s_p$ is zero, where $\pi_p$ is the projection $H^2(F, E[p]^{\otimes 2}) \rightarrow H^2(F, \mu_p)$. Thus the map we are interested in, is in effect
$$
T_F(A)/p \, \rightarrow \, H^2(F, {\rm sym}^2(E[p])).
$$
We will (without loss) use $s_p$, or $s_{p, F}$ to stress the dependence on $F$, to denote this map as well.

\bigskip

Thanks are due to some brief but helpful discussions over the years with Bloch, Colliot-Th\'el\`ene and Deligne, and to Raskind for a pertinent remark. The influence of Serre's works should be evident. This work was supported by a grant from the Simons Foundation. Finally, I thank the referee for a careful reading of the article, and for pointing out some relevant recent references and typos.

\bigskip

\section{A Local to Global result}

\medskip

For any number field
$F$ and any Gal$(\overline F/F)$-module $M$, set
$$
{\Sh}^i(F,M) : = \, {\text ker}(H^i(F,M) \rightarrow \prod_v
H^i(F_v, M)),
$$
which is known to be finite if $M$ is finite (\cite{Mi}, Theorem
4.10, page 70).

\medskip

\noindent{\bf Theorem C} \quad \it Let $E$ be an
elliptic curve over a number field $F$, and $p$ a prime. Then
\Sh$^2(F, {\rm sym}^2 E[p])$ is zero.
\rm

\medskip

When the image $G_p$ of Gal$(\overline F/F)$ in Aut$(E[p]) \simeq$
GL$(2, \F_p)$ is solvable, Theorem C is well known. We reduce the general non-solvable
case to when $G_p$ is SL$(2, \F_p)$ by making use of the knowledge
of possible subgroups of GL$(2, \F_p)$, and then appeal to a
result of J.P.~Serre relating the
cohomology of SL$(2,\F_p)$ to that of its Borel subgroup. See
below.

\medskip

\subsection{Preliminaries} \label{subsec1.1}

\medskip

For every finite Galois module $M$, put
$$
M^\vee :=\rm{Hom}(M, \overline F^\ast),
\leqno{(1.1.1)}
$$
equipped with the natural ${\ScrG}_F$--action. By the global {\it Tate
duality} (cf. Theorem 4.10, chap. I of \cite{Mi}), we know that
$\Sh^2 (F, M)$ is in duality with $\Sh^1 (F,
M^\vee)$. Moreover, if $\rho_p$ denotes the mod $p$ representation
of ${\ScrG}_F$ attached to $E$, the Cartier dual of sym$^2(\rho_p)$ identifies with Ad$(\rho_p)$,
which is sym$^2(\rho_p)\otimes\omega_p^{-1}$ (see \cite{MR} for example), where $\omega_p$ is the mod $p$ cyclotomic character.

Consequently it suffices to
prove the assertion that for any prime $p$,
$$
\Sh^1(F, {\rm Ad}(\rho_p)) \, = \, 0.\leqno(1.1.2)
$$

Let $K$ denote the smallest Galois extension of $F$ over which
$\rho_p$ becomes trivial (under the action of
${\ScrG}_{K}$). Let $G$ denote the Galois group of $K$ over
$F$. Then we have the following commutative diagram with exact
rows:
\small
$$
\begin{array}{ccccccc}
H^1 (G , {\rm Ad}(\rho_p))&\overset{\inf}{\hookrightarrow}&H^1 (F, {\rm Ad}(\rho_p))&\overset{\res}{\rightarrow}&H^1 (K, {\rm Ad}(\rho_p))\\
  \downarrow                 &                            & \downarrow         &          &\downarrow\\
\prod\limits_{v} \prod\limits_{w|v} H^1 (K_{w}/F_v ,
{\rm Ad}(\rho_p))&\overset{\inf}{\hookrightarrow}&\prod\limits_v H^1 (F_v , {\rm Ad}(\rho_p))&\overset{\res}{\rightarrow} &\prod\limits_v \prod\limits_{w|v}H^1
(K_{w}, {\rm Ad}(\rho_p))
\end{array}
\leqno{(1.1.3)}
$$
\normalsize
where \emph{inf} (resp. \emph{res}) denotes the inflation (resp.
restriction) map, $v$ runs over the places of $F$ and $w$ over
those of $K$, and the vertical maps are the global--to--local
maps defined by restriction. Since ${\rm Ad}(\rho_p)$ is acted upon
trivially by ${\ScrG}_{K}$, we have
$$
H^1 (K , {\rm Ad}(\rho_p)) \simeq \Hom ({\ScrG}_{K} ,
\mathbb{Z}/p ) \otimes_{\mathbb{Z}/p} {\rm Ad}(\rho_p)
\leqno{(1.1.4)}
$$
Similarly, $H^1 (K_{w}, {\rm Ad}(\rho_p)) \, \simeq \, \Hom
({\ScrG}_{K_{w}} , \mathbb{Z}/p )
\otimes_{\mathbb{Z}/p} {\rm Ad}(\rho_p)$. Then the right hand
vertical map is given by $j \otimes 1$, where
$$
j : \Hom ({\ScrG}_{K}, \mathbb{Z}/p)\, \to \, \prod_w \Hom
({\ScrG}_{w}, \mathbb{Z}/p ), \leqno{(1.1.5)}
$$
which is evidently injective. Thus in order to establish the assertion
1.1.2, it suffices to prove the following

\medskip

\noindent{\bf Proposition 1.1.6} \, \it The group
$$
\Sh^1_X (G , {\rm Ad}(\rho_p)):=Ker \left(H^1 (G , {\rm Ad}(\rho_p))
\overset{res}{\longrightarrow} \prod_{G^\prime\ \in\ X} H^1 (G^\prime, {\rm Ad}(\rho_p))\right)
$$
is zero, where $X$ denotes the set of decomposition groups $Gal
(K_{w}/F_v )$ in $G$. \rm

\medskip

\subsection{Hasse principle for subgroups of $GL_2
(\mathbb{F}_p )$} \label{subsec8.4}

\medskip

Now Proposition 1.1.6, (1.1.2) and hence Theoem C, will follow if we
establish the Hasse principle for any subgroup of $GL_2
(\mathbb{F}_p ).$ Following \cite{Mi}, p. 143-144, we will set, for any finite group $G$ and $G$-module $M$,
$$
H^1_{\ast}(G, M) \, = \, {\rm Ker}\left(H^1(G, M) \rightarrow \prod_C H^1(C, M)\right),
$$
where the product is over cyclic subgroups $C$ of $G$.

\medskip

\noindent{\bf Proposition 1.2.1} \, \it Let G be a subgroup of
$GL_2 (\F_p ).$ Then for any $p$-primary G-module M, we have
$$
H^1_\ast (G, M) =0.
$$
\rm

\medskip

\noindent{\bf Proof}. \, Let $B=U \rtimes T$ be the Borel subgroup
of $GL_2 (\mathbb{F}_p)$, with

$$
T= \left\{ \begin{pmatrix}a&0\\
                       0&b
\end{pmatrix} \biggm| a, b\in \mathbb{F}^\ast_p\right\}
$$
and
$$
U=\left\{ \begin{pmatrix}1&x\\
                  0&1
\end{pmatrix} \bigm| x\in \mathbb{F}_p\right\}.$$
\medskip

Suppose G has order prime to $p$. Then G will be a
subgroup of a standard maximal torus, which will be either the split torus $T$ above, or the
anisotropic torus defined by the multiplicative group of the quadratic extension
$\F_{p^2}$ of $\F_p$. The assertion  is clear in
this case.

So we may assume from here on that G has order divisible by
$p$. We can then appeal to the following result of
J.~P.~Serre:

\medskip

\noindent{\bf Proposition 1.2.2} (\cite{Serre1}, sec. 2.4, Prop.
15) \it Let $G$ be a subgroup of $GL_2 (\mathbb{F}_p )$ with
$p \bigm| |G|$. Then exactly one of the following must happen:

(i) $G$ is conjugate to a subgroup of the Borel subgroup $B$

(ii) G contains $SL_2 (\mathbb{F}_p )$. \rm

\medskip

If we are in case (i), G is solvable and the theorem is known cf. \cite{Mi}, chap. 1,
Prop. 9.5 (b).

So we may assume that
$$ G \supset SL_2 (\F_p ).
\leqno{(1.2.3)}
$$

Then we have a short exact sequence
$$
1 \to SL_2 (\F_p) \to G \overset{\rm det}{\rightarrow} \overline{G} \to 1
$$
with
$$
\overline{G} \subset \F^\ast_p.
$$

This gives rise to the following diagram with exact rows
\small
$$
\begin{CD}
0 \to H^1 (\overline{G}, M^{SL_2(\F_p)}) @>inf>> H^1 (G, M) @>\rm{res}>> H^1 \left(SL_2(\F_p), M\right)\\
               @VVV                             @VVV                                             @VVV \\
0 \to \prod\limits_{C_1} H^1 (C_1 , M^{SL_2(\F_p)}) @>inf>>
\prod\limits_C H^1 (C, M) @>\rm{res}>> \prod\limits_C H^1 (C \cap
SL_2(\F_p), M)
\end{CD}
\leqno{(1.2.4)}
$$
\normalsize
where $C$ runs over the cyclic subgroups of $G$, and $C_1$ denotes
the image of $C$ in $\overline{G}$. We note that any cyclic
subgroup of $SL_2(\F_p)$ (resp. $\overline{G}$) is of the form
$C\cap G$ (resp. $C_1$).

Moreover, since $\overline{G}$ is cyclic of order prime to $p$,
$H_*^1 (\overline{G}, M^{SL_2(\F_p)})$ is zero. Thus in order
to prove Proposition 1.2.1, it suffices to establish the following

\medskip

\noindent{\bf Proposition 1.2.5} \it For any $p$-primary
$SL_2(\F_p)$-module M, we have:
$$
H_*^1 (SL_2(\F_p), M)=0.
$$
\rm

\medskip

\noindent{\bf Proof}. \, Let B also denote, by abuse of
notation, the Borel subgroup of $SL_2(\F_p)$. Now we will
appeal to another result of J.--P.~Serre: \rm

\medskip

\noindent{\bf Theorem 1.2.6} (\cite{Serre3}, Proposition 13
and Corollary, page 555) \, \it For any $p$-primary
$SL_2(\F_p)$-module M, the restriction map
$$
H^* (SL_2(\F_p), M)\overset{\alpha}{\rightarrow} H^* (B, M)
$$
is injective. \rm
\medskip

The injectivity of $\alpha$ reduces the problem to showing that
$$
H_*^1 (B, M)=0 \leqno{(1.2.7)}
$$
for any $p$-primary module $M$. Since B is solvable, we are done. \qed
\medskip

This concludes the proof of Proposition 1.2.1, Proposition 1.1.6 and
Theorem C.

\qed

\bigskip

\section{Proof of Theorem A}

\medskip

Let $E/F, p$ as in Theorem A. Pick any $\theta \in T_F(A)/p$. In view of Theorem C (of the section above), it suffices to show for {\it every} place $v$ of $F$, the image (under restriction) of $s_{p,F}(\theta)$ in $H^2(F_v, {\rm sym}^2(E[p]))$ is zero. By the naturality of the Galois symbol,
this image is the same as $s_{p, F_v}(\theta_v)$, where $\theta_v$ is the image of $\theta$ in $T_{F_v}(A)$.

First consider when $v$ archimedean. If it is complex, the target group $H^2(F_v, {\rm sym}^2(E[p]))$ is zero, so there is nothing to prove. When $v$ is real, since the Galois group is $\Z/2$, the image is zero when $p$ is odd.

Let us take $p$ from here on to be odd and a prime where $E$ has good ordinary reduction. Furthermore, let us also assume that $p$ does not divide the discriminant of $F$. All this only eliminates a finite number of primes.

Now let $v$ be a finite place with residue characteristic prime to $p$. Then (as seen in \cite{MR}), $s_{p, F_v}(\theta_v)$ lands in $H^2(F_v^{\rm nr}/F, {\rm sym}(E[p])$, which is zero. (Here $F_v^{\rm nr}$ denotes as usual the maximal unramified extension of $F_v$, so that Gal$(F_v^{\rm nr}/F) \simeq \hat\Z$.)

It remains to consider the main case when $p$ is the residual characteristic of $F_v$. In the case, by Theorem A, parts (a), (b) of \cite{MR}, the image of $s_{p, F_v}$ is at most one-dimensional, and moreover, this dimension is $1$ iff we have

\medskip

\noindent($\ast_v$) \, $\mu_p \subset F_v$ {\it and} $F_v(E[p])$ is unramified over $F_v$ with the prime-to-$p$ part of $[F_v(E[p]):F_v]$ being $\leq 2$.

\medskip

As $[F_v(\mu_p):\Q_p]$ is $\geq p-1$, a necessary condition for $(\ast_v)$ to occur is for the Galois closure $\tilde F$ of $F$ over $\Q$ to have degree $\geq p-1$ (over $\Q$). Thus $s_{p, F_v}(\theta_v)$ is zero if $p$ chosen above is $> [\tilde F : F]$.

Consequently, $s_{p,F}(\theta)$ is zero for all but finitely many primes $p$ where $E$ has good ordinary reduction.

\qed

\bigskip

\section{Special classes in the albanese kernel via $p$-adic approximation}

\medskip

One sees via archimedean considerations that if
we have an elliptic curve $E$ over a number field $F$ with a real
embedding such that $E[2]$ is rational over $F$, then for a
suitable $2$-torsion point $P$, the class $\langle P, P\rangle$ is
non-trivial in $T_F(A)/2$. So we will restrict ourselves in this
chapter to the construction of non-trivial classes in $T_K(A)/p$
for suitable extensions $K/F$ when $p$ is an {\it odd} prime.
Since $T_F(A)$ is expected to be a finite group, which is
compatible with our theorem in the previous chapter, for arbitrary
$p$ we will need to go to finite extensions $K/F$ to get
non-trivial cycle classes in the Albanese kernel.
We acknowledge here with thanks some helpful remarks of Deligne regarding this matter some years back.

\medskip

The following result is a more precise version of Theorem B of the
Introduction.

\noindent{\bf{Theorem D}} \, \it Let $F$ be a number field.
\begin{enumerate}
\item[a)] \, Let $E$ be an elliptic curve over $F$. Let $p$ be
an odd prime such that $E$ has good ordinary reduction at some place
$v$ of $F$ dividing $p$. Then there exists a finite extension $K$ of $F$
containing $E[p]$, and points $P, Q\ \epsilon\ E(K)$, with $P \in
E[p]$, such that
$$
\langle P,Q\rangle \, \neq \, 0 \, \in T_K (E\times E)/p.
$$
In fact we can take $K$ such that $[K:F(E[p])]\leq 4$. Moreover,
after possibly replacing $K$ by a finite extension unramified
above $v$, we may choose $P$ to be a $p$-power torsion point.
\item[b)] \,  Let $p$ be an odd
prime. Then there exists a number field $K$ containing $F$,
an elliptic curve E defined over F, and points $P, Q\ \epsilon\
E(K)$ such that $\langle P, Q\rangle \neq 0$ in $T_K(A)/p.$
\end{enumerate}
\rm

\medskip

This Theorem answers a question of J.-L.Colliot-Th\'el\`ene (private communication from some years back). All the classes $\langle P, Q\rangle$ should be torsion, but we do not
know how to show this, except of course when we start with a
torsion point $Q$. One can bound the degree of the finite extension over which
we may choose $P$ to be a $p$-power torsion point.

\medskip
When for CM elliptic curves $E$, non-trivial global $0$-cycles on the albanese kernel of $E \times E$ have already been constructed by a very different method \cite{GK}, Theorem 1.7, which actually exist in a smaller degree extension; see also the references therein. (I am thankful to the referee for pointing this out.)

\medskip

\subsection{A lemma from non-archimedean geometry}
\label{subsec10.2}

\medskip

We need the following fact:

\medskip

\noindent{\bf Lemma 3.1.1} \, \it Let $f:X \rightarrow Y$, be a
morphism of smooth algebraic varieties X, Y defined over a local
field F. Let $P \epsilon X(F), P' \epsilon Y(F)$ and $P'=f(P).$
Assume that $f$ is $\acute{e}$tale (in the sense of algebraic
geometry) at P. Then there exist, in the analytic topology,
neighborhoods $U(P) \subset X(F)$ of P and $U'(P') \subset Y(F)$
of $P'$ such that $f:U(P)\rightarrow U'(P')$ is an analytic
isomorphism.

Moreover, if Y is the affine line, then we can take $U'(P')$ to be
an open unit ball around $P'$.

\rm

\medskip

For a proof, see \cite{BGR}, page 302, Prop. 5.

\medskip

\subsection{Approximation} \label{subsec10.3}

We will use the fact that $E$ can be realized as a double cover of
$\mathbb P^1$ with $E[2]$ being the inverse image of the
ramification locus.

\medskip
\noindent \textbf{$Proposition$ 3.2.1} Let $P_1 \epsilon E(F_v )$,
$P_1 \notin E[2]$. Then there exists a finite extension $K\supset
F$, a place w of K such that $K_w = F_v$, and a point $P\
\epsilon\ E(K)$ such that $P \equiv P_1 \left(\mod\, p
E(K_w)\right)$. Moreover, if $F\supset E[p],$ we can take K such
that $[K:F]\leq 2$.

\medskip

{\it Proof}. \, First of all we construct an analytic neighborhood
$\tilde{U}(P_1) \subset E(F_v)$ of $P_1$ such that for every
$P^\ast \epsilon \widetilde{U}(P_1)$ we have $P^\ast-P_1\
\epsilon\ p E(F_v)$. For this we apply Lemma 1.2.1 to the
everywhere \'etale morphism $p:E\to E$ (multiplication by p) and
to the points $R \in E[p]$ and $0 \in E$ with $R \neq 0$. (Note
that since $E[p] \subset F$ we have such a point $R \in E(F)
\subset E(F_v)$). By Lemma 3.1.1 we have neighborhoods
$\widetilde{U}(R) \subset E(F_v)$ of R and
$\widetilde{U}(0)\subset E(F_v)$ of 0 such that
$p:\widetilde{U}(R)\overset{\sim}{\longrightarrow} \widetilde{U}(0)$. Now take
$\widetilde{U}(P_1)=\widetilde{U}(0)+P_1$. This satisfies our
requirement.

Next, without loss of generality we can assume that our curve is
given by a Weierstrass equation $Y^2 =\phi(X) \in F(X)$ with
$\phi(X)$ of degree 3. Let $P_1 =(x_1 , y_1 )$. Consider the
morphism $\pi : E\to \mathbb{P}^2$ given by $\pi (x, y)=x$. Since
$P_1 \notin E[2]$ the map $\pi$ is \'{e}tale in $P_1$. Apply
$Lemma$ 3.1.1 to $\pi$ and to the points  $P_1 =(x_1 , y_1 )$ on
$E$ and $x_1 \in\mathbb{P}^1$. By this $Lemma$ there exist
analytic neighborhoods $U(P_1 )$ of $P_1$ in $E(F_v )$ and $U(x_1
)$ of $x_1$ in $\mathbb{P}^1 (F_v )$ such that $\pi :U(P_1 )
\overset{\sim}{\longrightarrow} U(x_1 )$. Moreover, we can clearly assume that
$U(P_1 ) \subset \widetilde{U}(P_1 )$, the analytic neighborhood
constructed above. (Also we can take for $U(x_1 )$ the unit ball.)

Now take $x\ \epsilon\ F$ such that $x\ \epsilon\ U(x_1 )$. The
equation $Y^2 =\phi(x)$ has two roots $y'$ and $y''$, which are
distinct because $\pi$ is \'{e}tale above $U(x_1 )$. One of the
roots, say $y'$, is such that $(x, y')$ lies in $U(P_1 ) \subset
E(F_v)$.Write $y = y'$ and let $K=F(y)$, so that $[K:F] \leq 2$.
Since $y=y'$ is in $F_v$ we have $K\subset F_v$, and hence there
is a place $w$ of $K$ extending $v$ such that $K_w =F_v$. Now the
point $P=(x, y) \in E(K)$ is in $U(P_1 )\subset \widetilde{U}(P_1
)$, so $P \equiv P_1\big( (\mod\ p E(K_w)\big)$ and $P$ clearly
fulfills the requirements. This proves the $Proposition$.

\qed

\medskip

\subsection{Proof of Theorem D} \label{subsec10.4}

\medskip

First we prove part a). Here $E$ is an elliptic curve over a
number field $F$, $p$ an odd prime, and $v$ a place of $F$ above
$p$ such that $E$ has good ordinary reduction at $v$. Replacing
$F$ by $F(E[p])$ if necessary, we can assume for establishing part
a) of Theorem D, that $F\supset E[p]$. Then $\mu_p$ is also in $F
\subset F_v$, and so by part (c) of \cite{MR}, Theorem A,
we can find $P_1, Q_1, \in E(F_v)$, not in $E[2]$
(because $p \ne 2$), such that $\langle P_1, Q_1\rangle \neq 0$ in
$T_{F_v} (E\times E)/p$. Over a finite extension we can also
choose $P_1$ to be a $p$-power torsion point, if so desired.

Now apply Proposition 3.2.1 to both $P_1$ and $Q_1$ and get a
field $K\supset F$ with $[K:F]\leqq 4$ a place $w$ of K dividing
$v$ such that $K_w = F_v$, together with points $P, Q \in E(K)$
such that $P \equiv P_1, \rm{resp}.\ Q\equiv Q_1 \left( \mod\ p
E(K_w )\right)$. Then we have $\langle P, Q\rangle = \langle P_1 ,
Q_1 \rangle$ in $T_{K_w}(E\times E)/p$. Since $T_{K_w}(E\times
E)=T_{F_v}(E\times E)$ and $\langle P_1 , Q_1 \rangle \neq 0$ in
$T_{K_w}(E\times E)/p$.

\medskip

Now on to part b). By hypothesis, p is a prime number $ \neq 2$.
There exists an elliptic curve $E$ defined over $\mathbb{Q}$ which
has ordinary reduction at $p$. In fact there exist infinitely many
such curves, and we can even take E to be CM. Put
$F=\mathbb{Q}(E[p])$.

By Theorem A of \cite{MR}, $\exists P_1 , Q_1, \in E(F_v )$ such that
$\langle P_1 , Q_1 \rangle \neq 0$ in $T_{F_v}(E\times E)/p$.
Since $E[p] \subset F$ we take $P=P_1$. Now applying $Proposition$
3.2.1 there exists an extension $K\supset F$ with $[K:F]\leq 2$,
a place $w$ of $K$ such that $w$ divides $v$ and $K_w = F_v$, and
a point $Q \in E(K)$ such that $Q\equiv Q_1 \left(\mod\ p E(F_v
)\right)$. Then $\langle P, Q \rangle = \langle P_1 , Q_1 \rangle$
in $T_{F_v}(E\times E)/p$, and hence $\langle P, Q\rangle \neq 0$
in $T_{F_v}(E\times E)/p$. This proves part b) (with the field $K$
and $E$ over $K$).

\qed

\noindent Dinakar Ramakrishnan

\noindent Department of Mathematics

\noindent California Institute of Technology

\noindent Pasadena, CA 91125, USA

\noindent \texttt{dinakar@caltech.edu}

\end{document}